\documentclass{amsart}

\usepackage{cite}
\usepackage{mathrsfs}
\usepackage{subfig}
\usepackage{graphicx}
\usepackage{multirow}
\usepackage{url}
\usepackage{microtype}
\usepackage{hyperref}

\def\al{\alpha} 
\def\be{\beta}

\def\th{\theta}

\def\la{\lambda}

\let\on=\operatorname

\def\R{\mathbb{R}}

\newtheorem{theorem}{Theorem}

\theoremstyle{definition}
\newtheorem{definition}{Definition}
\newtheorem{remark}{Remark}

\begin{document}

\title{Curve Matching with Applications in Medical Imaging}
\author{Martin Bauer}
\address{Fakult\"at f\"ur Mathematik, Universit\"at Wien}
\author{Martins Bruveris}
\address{Department of Mathematics, Brunel University London}
\author{Philipp Harms}
\address{Department of Mathematics, ETH Z\"urich}
\author{Jakob M\o ller-Andersen}
\address{Department of Applied Mathematics and Computer Science, Technical University of Denmark}
\thanks{All authors were partially supported by the Erwin Schr\"{o}dinger Institute programme: Infinite-Dimensional Riemannian Geometry with Applications to Image Matching and Shape Analysis. Martin Bauer was supported by FWF-grant no. P24625.}
\date{June 2015}
\subjclass{Primary 58B20; Secondary 62H25, 62H30}
\keywords{Curve matching, Sobolev metrics, Riemannian shape analysis, discrete geodesics, minimizing geodesics}

\begin{abstract}
In the recent years, Riemannian shape analysis of curves and surfaces has found several applications in medical image analysis. In this paper we present a numerical discretization of second order Sobolev metrics on the space of regular curves in Euclidean space. This class of metrics has several desirable mathematical properties. We propose numerical solutions for the initial and boundary value problems of finding geodesics. These two methods are combined in a Riemannian gradient-based optimization scheme to compute the Karcher mean. We apply this to a study of the shape variation in HeLa cell nuclei and cycles of cardiac deformations, by computing means and principal modes of variations.
\end{abstract}

\maketitle

\section{Introduction}

The comparison and analysis of geometric shapes plays an important role in medical imaging \cite{Younes2012}, biology \cite{Dryden1998} as well as many other fields \cite{Krim2006}. Spaces of geometric shapes are inherently nonlinear. To make standard methods of statistical analysis applicable, one can linearize the space locally around each shape. This can be achieved by introducing a Riemannian structure, which is able to describe jointly the global nonlinearity of the space as well as its local linearity. 

Over the past decade, Riemannian shape analysis has become an active area of research in pure and applied mathematics. Driven by applications, a variety of spaces, equipped with different Riemannian metrics, have been used to represent geometrical shapes and their attributes. Ideally, the particular choice of metric should be dictated by the data at hand rather than by mathematical or numerical convenience. This leads to the task of developing efficient numerical methods for the statistical analysis of shapes for general and flexible classes of metrics. 

The topic of this paper are second order Sobolev metrics on the space of regular, planar curves. This space and its quotients by translations, rotations, scalings, and reparametrizations are important and widely used spaces in shape analysis. Second order Sobolev metrics are mathematically well-behaved: the geodesic equation is globally well-posed, any two curves in the same connected component can be connected by a minimizing geodesic, the metric completion consists of all $H^2$-immersions, and the metric extends to a strong Riemannian metric on the metric completion \cite{Bruveris2014, Bruveris2014b_preprint}. 

We provide the first numerical implementation of the initial and boundary value problems for geodesics with respect to these metrics and apply them to medical data.\footnote{Implementations of second order Sobolev metrics can be found in \cite{Vialard2014_preprint} and \cite{Bauer2014c}, but none of the previous implementations seem to have been used in applications.} Our implementation allows us to compute Karcher means, principal components, and clusters of curves in reasonable time. The parameters in the metric can be chosen freely because we do not rely on transforms which exist only for special choices of parameters \cite{Jermyn2011, Michor2008a}. We are also able to factor out the action of the finite-dimensional translation and rotation groups.

We illustrate the behaviour of the metrics in two medical applications. In Sect.~\ref{hela} we use images of HeLa cell nuclei from \cite{Boland2001} and compute the mean shape of the nucleus and the principal modes of shape variation; in Sect.~\ref{cardiac} we use traces of cardiac image sequences from \cite{tobongomez, Mcleod} and compare the mean shape of patients with Tetralogy of Fallot with the mean shape of a control group.

The paper is structured as follows: Sect.~\ref{math} contains the definition of second order Sobolev metrics and some relevant mathematical background. Section \ref{disc} describes our discretization of the geodesic equation and the Riemannian energy functional. Finally, in Sect.~\ref{app}, we apply the metrics to medical imaging data. 

\section{Mathematical Background}
\label{math}

In this article we center our attention on the space of smooth, regular curves with values in $\R^d$,
\begin{align}
\on{Imm}(S^1,\mathbb R^d)=\left\{c\in C^{\infty}(S^1,\mathbb R^d)\colon \forall \th \in S^1, c_\th(\th) \neq 0 \right\}\enspace, 
\end{align}
where $\on{Imm}$ stands for \emph{immersion}. As an open subset of the Fr\'echet space $C^\infty(S^1,\R^d)$, it is a Fr\'echet manifold. Its tangent space, $T_c\on{Imm}(S^1,\R^d)$, at any curve $c$ is the vector space $C^\infty(S^1,\R^d)$. 

We denote the Euclidean inner product on $\mathbb R^d$ by $\langle\cdot,\cdot\rangle$. Moreover, for any fixed curve $c$, we denote differentiation and integration with respect to arc length by $D_s=\frac{1}{|c_\theta|}\partial_{\theta}$ and $ds=|c_\theta|d\theta$ respectively. 

\begin{definition}
\label{def:sobolev_metric}
A second order Sobolev metric on $\on{Imm}(S^1,\mathbb R^d)$ is given by
\begin{equation}\label{def:met}
G_c(h,k) = \int_{S^1} a_0\langle h,k \rangle+a_1\langle D_s h,D_s k \rangle+a_2\langle D_s^2 h,D_s^2 k \rangle \,ds \enspace,
\end{equation}
where $a_0, a_2 > 0$, $a_1 \geq 0$, and $h,k \in T_c\on{Imm}(S^1,\mathbb R^d)$ are tangent vectors.
\end{definition}

\begin{remark}[Invariance of the metric]
The metric $G$ is invariant with respect to translations, rotations, and reparametrizations by diffeomorphisms of $S^1$, but not with respect to scalings. It can be made scale-invariant by introducing weights that depend on the length $\ell_c$ of the curve $c$. A scale-invariant metric is given by
\begin{equation}
\widetilde G_c(h,k) = \int_{S^1} \frac{a_0}{\ell^3_c}\langle h,k \rangle+\frac{a_1}{\ell_c} \langle D_s h,D_s k \rangle+a_2\ell_c \langle D_s^2 h,D_s^2 k \rangle\, ds\enspace.
\end{equation}
Depending on the application, it can be desirable to factor out some of these isometry groups. For example, we consider the HeLa cell nuclei in Sect.~\ref{hela} as curves modulo translations and rotations. For the traces of cardiac images in Sect.~\ref{cardiac}, however, both the absolute position of the curve and its orientation have intrinsic meaning and should therefore not be factored out.
\end{remark}

\begin{figure}
	\centering
	\includegraphics[width=0.95\textwidth]{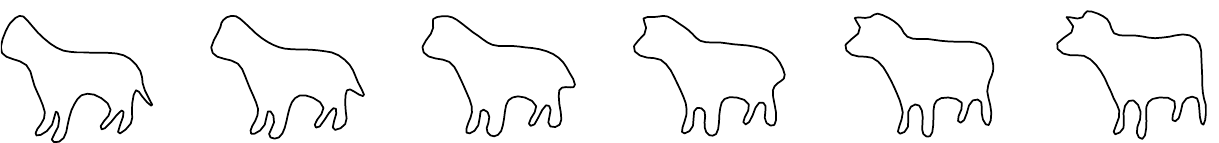}
	\includegraphics[width=0.95\textwidth]{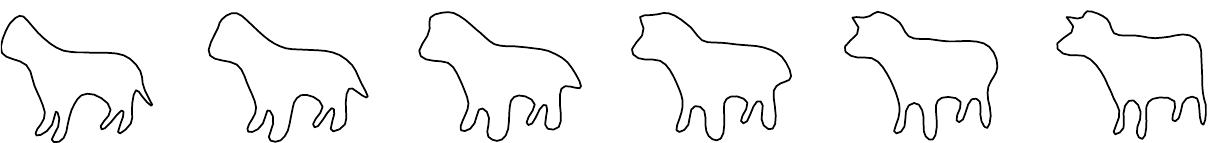}
	\includegraphics[width=0.95\textwidth]{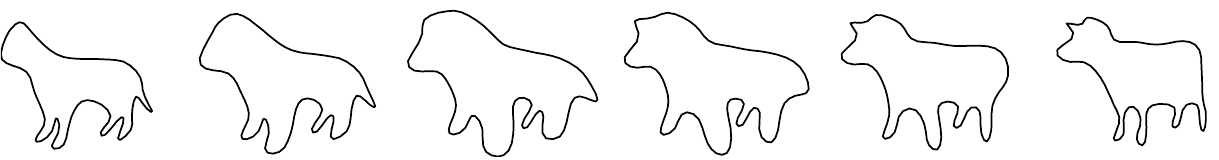}
	\caption{Geodesics between a cat and a cow. The metric parameter $a_2$ in the first row is increased by a factor 10 in the second and a factor 100 in the third row.}
	\label{fig:cat2cow_constants}
\end{figure}

\begin{figure}
\centering
	\includegraphics[width=.15\textwidth]{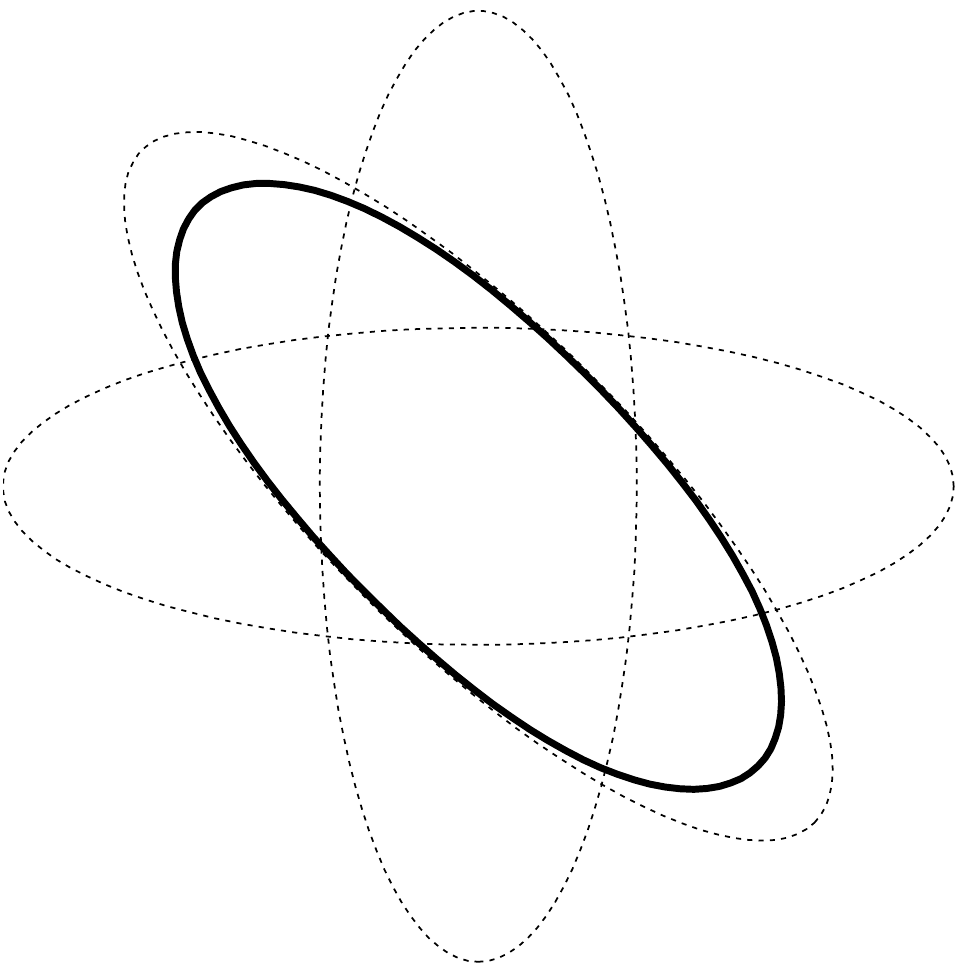}
	\hspace{1cm}
	\includegraphics[width=.15\textwidth]{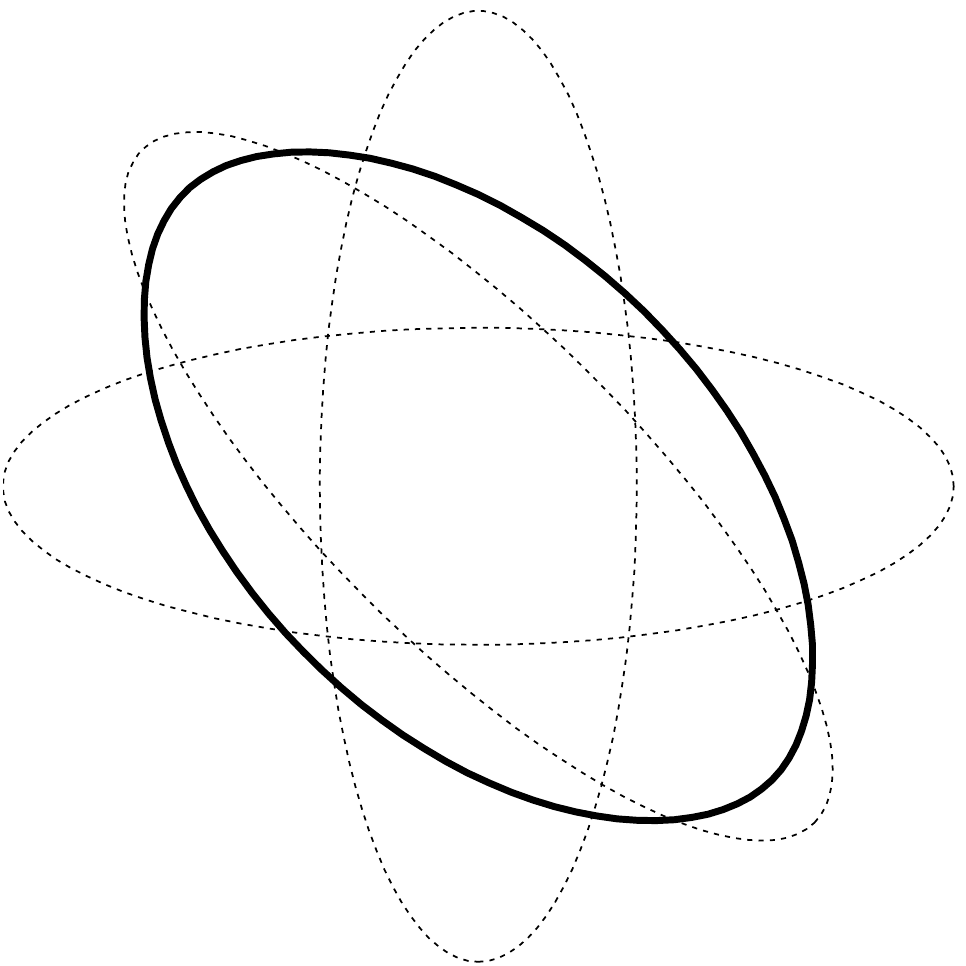}
	\hspace{1cm}
	\includegraphics[width=.15\textwidth]{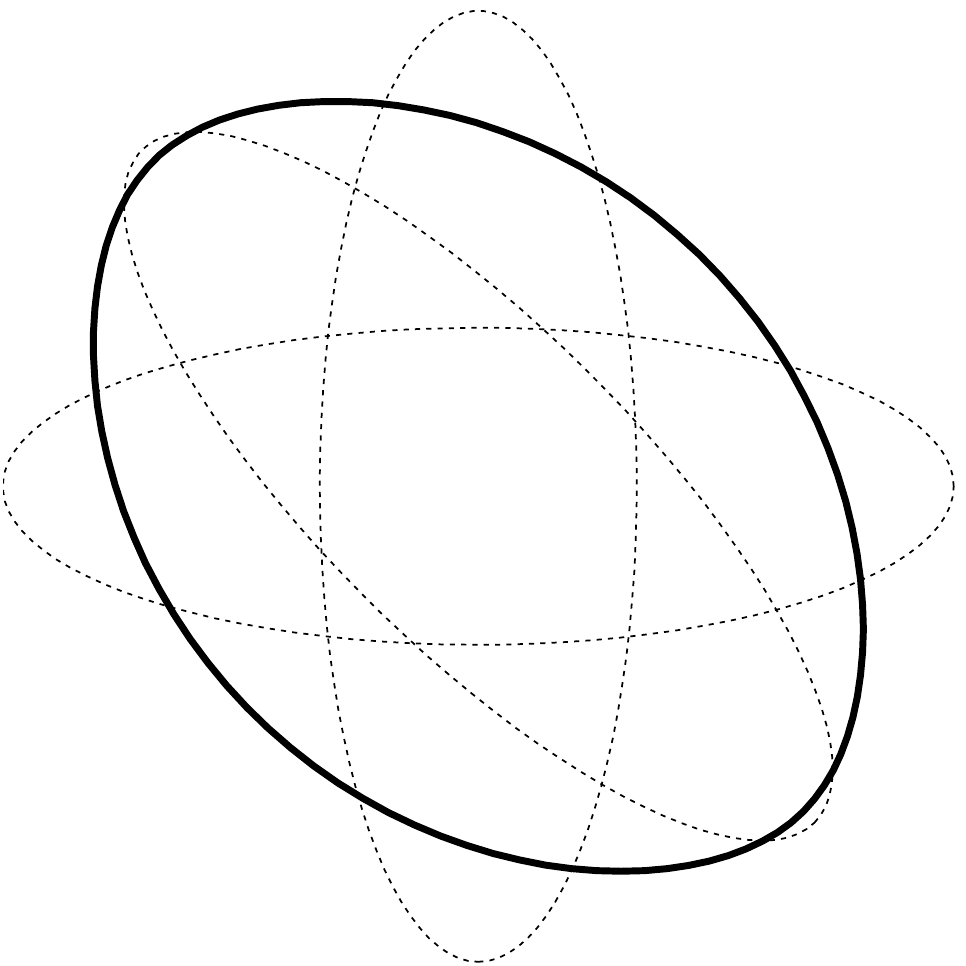}
	\caption{Karcher means (bold) of three rotated ellipses (dashed). The parameter $a_2$ of the metric used in the first figure is increased by a factor 10 in the second and by a factor 100 in the third figure.}
	\label{fig:Karcher_constants}
\end{figure}

\begin{remark}[Choosing the constants]
The constants $a_0, a_1$ and $a_2$ in the definition of the metric determine the relative weight of the $L^2$-, $H^1$-, and $H^2$-parts. Their choice is a non-trivial and important task, and it should, in each application, be informed by the data at hand. The influence of different parameter values on geodesics and Karcher means can be seen in Fig.~\ref{fig:cat2cow_constants} and Fig.~\ref{fig:Karcher_constants}.
\end{remark}

\begin{remark}[Generalizations to higher-dimensions]
Sobolev metrics have natural generalizations to manifold-valued curves, embedded surfaces, and more generally spaces of immersions of a compact manifold $M$ into a Riemannian manifold $N$; see \cite{Bauer2011a, Bauer2011b} for details and \cite{Bauer2014} for a general overview.
\end{remark}

Deformations of curves can be seen as smooth paths $c\colon\![0,1]\to\on{Imm}(S^1,\mathbb R^d)$. Their velocity is $c_t$, the subscript $t$ denoting differentiation.

\begin{definition}[Geodesic distance]
The length of a path of curves $c$ is 
\begin{align}
L(c) = \int_0^1 \sqrt{G_{c(t)}(c_t(t),c_t(t))} \, dt\enspace.
\end{align}
The distance between two curves in $\on{Imm}(S^1,\R^d)$ (with respect to the metric $G$) is the infimum of the lengths of all paths connecting these curves. 
\[
\on{dist}(c_0, c_1) = \inf_{\substack{c(0) = c_0, c(1) = c_1}} L(c)\,.
\]
Geodesics are locally distance-minimizing paths.
\end{definition}

Geodesics can be described by a partial differential equation, called the geodesic equation. It is the first order condition for minima of the energy \begin{equation}
E(c) = \frac12\int_0^1 G_{c(t)}\big(c_t(t),c_t(t)\big) \, dt\enspace.
 \label{eq: EnergyFunctional}
\end{equation}
Recently, some local and global existence results for geodesics of Sobolev metrics were shown in \cite{Michor2007,Bruveris2014,Bruveris2014b_preprint}. Since they provide the theoretical underpinnings for the numerical methods presented in this paper, we summarize them here.

\begin{theorem}[Geodesic equation]
Let $a_0, a_2 > 0$ and $a_1 \geq 0$. The geodesic equation of the metric $G$, written in terms of the momentum 
$p=|c'| (a_0 c_t-a_1 D^2_s c_t+a_2 D_s^4c_t)$, is given by
\begin{align}
\partial_t p = &-\frac{a_0}2 |c_\theta| D_s(\langle c_t,c_t \rangle D_sc) + \frac{a_1}{2}|c_\theta| D_s (\langle D_s c_t,D_s c_t \rangle D_sc) \nonumber
\\&
-\frac{a_2}{2}|c_\theta|D_s (\langle D^3_s c_t,D_s c_t \rangle D_sc)+\frac{a_2}{2}|c_\theta|D_s (\langle D^2_s c_t,D^2_s c_t \rangle D_sc)\,.
\label{eq: GeodesicEquation}
\end{align}
Given any initial condition $(c_0, u_0) \in T\on{Imm}(S^1,\R^2)$, the solution of the geodesic equation exists for all time.

If, however, $a_0, a_1 > 0$ and $a_2=0$, then $G$ is a first order Sobolev metric. Its geodesic equation is locally, but not globally, well-posed.
\end{theorem}

\begin{remark}[Comparison to elastic metrics]
Closely related are elastic metrics \cite{Mio2007}, which in the planar case are given by
\begin{equation}
\label{elastic}
 G_c(h,k) = \int_{S^1} a^2\langle D_s h,n \rangle \langle D_s k,n \rangle+b^2\langle D_s h,v \rangle\langle D_s k,v \rangle ds\enspace.
\end{equation}
Here $a,b$ are constants and $v,n$ denote the unit tangent and normal vectors to $c$. Two special cases deserve to be highlighted: for $a=1$, $b=\frac 12$  \cite{Jermyn2011} and $a=b$ \cite{Michor2008a} there exist nonlinear transforms, the square root velocity transform and the basic mapping, that greatly simplify numerics. Both of these metrics have been applied to a variety of problems in shape analysis.

We note that the elastic metric with $a=b$ corresponds to a first order Sobolev metric as in Def.~\ref{def:sobolev_metric} with $a_0=a_2=0$ and $a_1 = a^2=b^2$. As it has no $L^2$-part, it is a Riemannian metric only on the space of curves modulo translations.
\end{remark}

\section{Discretization}
\label{disc}

We propose to represent paths of curves as tensor product B-splines
\begin{equation}\label{eq: TensorProductPath}
c(t, \th) = \sum_{i=1}^{N_t} \sum_{j=1}^{N_\th} c_{i,j} B_i(t) C_j(\th)
\end{equation}
of degree $n_t$ in time ($t$) and $n_\th$ in space ($\theta$), respectively, using $N_t$ and $N_\th$ basis functions and controls $c_{i,j} \in \R^d$. The spatial basis functions $C_j$ are defined on the uniform knot sequence $\theta_j = \frac{j-n}{2\pi\cdot N_\theta}$, $0\leq j\leq 2 n+N_\theta$, and satisfy periodic boundary conditions. The time basis functions $B_i$ are defined using uniform knots on the interior of the interval $[0,1]$ and with full multiplicity on the boundary. Full multiplicity of the boundary knots ensures that fixing the initial and final curves is equivalent to fixing the control points $c_{1,j}$ and $c_{N_t,j}$. Regarding smoothness, we have $B_i \in C^{n_t-1}([0,1])$ and $C_j \in C^{n_\th-1}([0,2 \pi])$.

The proposed discretization scheme is independent of the special form of the metric, and can be applied equally to the scale-invariant versions as well as to the family of elastic metrics.

\subsection{Boundary Value Problem for Geodesics}
\label{sec: GeodesicBVP}
To evaluate the energy \eqref{eq: EnergyFunctional} on discretized paths \eqref{eq: TensorProductPath}, we use Gaussian quadrature on each interval between subsequent knots. The resulting discretized energy $E_{\on{discr}}(c_{1,1},\dots,c_{N_t,N_\th})$ is minimized over all discrete paths \eqref{eq: TensorProductPath} with fixed initial and final control points $c_{1,j}$ and $c_{N_t,j}$. This is a finite-dimensional, nonlinear, unconstrained\footnote{We neglect the condition $c_\th(t,\th) \neq 0$. All smooth paths violating this condition have infinite energy because of the completeness of the metric. By approximation, spline paths violating this condition have very high energy and will, in practice, be avoided during optimization.} optimization problem. 

We solved this problem using either Matlab's \texttt{fminunc} function with gradients computed by finite differences, or using AMPL \cite{fourer2002ampl} and the Ipopt solver \cite{wachter2006implementation} with gradients computed by automatic differentiation. By classical approximation results \cite{Schumaker2007}, the discrete energy of discrete paths converges to the continuous energy of smooth paths as the number of control points tends to infinity. Hence we expect to obtain close approximations of geodesics. Establishing rigorous convergence results for our discretization will be the subject of future work.

\subsection{Computing the Karcher Mean}
\label{karcher}
The Karcher mean $\overline{c}$ of a set $\{c_1,\dots,c_n\}$ of curves is the minimizer of
\begin{equation}
\label{eq: KarcherEnergy}
F(c) = \frac{1}{n} \sum_{j=1}^n \on{dist}(c, c_j)^2\enspace.
\end{equation}
It can be calculated by a gradient descent on $(\on{Imm}(S^1,\R^d), G)$. Letting $\on{Log}_{c}c_j$ denote the Riemannian logarithm, the gradient of $F$ with respect to $G$ is \cite{Pennec2006b}
\begin{equation}
\on{grad}^GF (c) = \frac 1n \sum_{j=1}^n \on{Log}_c c_j\enspace.
\end{equation}

\subsection{Initial Value Problem for Geodesics}
To calculate the Karcher mean by gradient descent, one has to repeatedly solve the geodesic equation \eqref{eq: GeodesicEquation}. To this aim, we use the time-discrete variational geodesic calculus \cite{Rumpf2014}. Given three curves $c_0,c_1,c_2$, one defines the discrete energy 
\begin{equation}
E_2(c_0,c_1,c_2) = G_{c_0}(c_1-c_0,c_1-c_0) + G_{c_1}(c_2-c_1,c_2-c_1)\enspace.
\end{equation}
A 3-tuple $(c_0,c_1,c_2)$ is a discrete geodesic if it is a minimizer of the discrete energy $E_2$ with fixed endpoints $c_0,c_2$. The discrete exponential map is defined as follows: $c_2 = \on{Exp}_{c_0} c_1$, if $(c_0,c_1,c_2)$ is a discrete geodesic, in other words, if $c_1 = \on{argmin} E_2(c_0,\cdot, c_2)$. To find $c_2$, we differentiate the discrete energy $E_2$ with respect to $c_1$ and solve the resulting system of nonlinear equations,
\begin{equation}
2G_{c_0}(c_1-c_0,\cdot) - 2 G_{c_1}(c_2-c_1,\cdot) + D_{c_1}G_\cdot(c_2-c_1,c_2-c_1)=0\enspace.
\end{equation}
We use the solver \texttt{fsolve} in Matlab to solve this system of equations. This procedure is repeated depending on the desired time-resolution.

While the convergence results in \cite{Rumpf2014} do not apply to our setting, we found good experimental agreement with the solutions of the boundary value problem.

\section{Applications}
\label{app}

\subsection{Hela Cells}
\label{hela}

In our first example we want to characterize nuclear shape variation in HeLa cells. We use fluorescence microscope images of HeLa cell nuclei\footnote{The dataset was downloaded from \url{http://murphylab.web.cmu.edu/data}.} (87 images in total). The acquisition of the cells is described in \cite{Boland2001}. A similar study on this dataset was performed using the LDDMM framework in \cite{Rohde2008, Rohde2008b}; applying intrinsically defined Sobolev metrics to the same problem will provide a complimentary point of view.

\begin{figure}
\centering
\begin{tabular}{cc}
\multirow{2}{*}{\smash{\raisebox{-0.5\height}{
\includegraphics[width=0.5\textwidth]{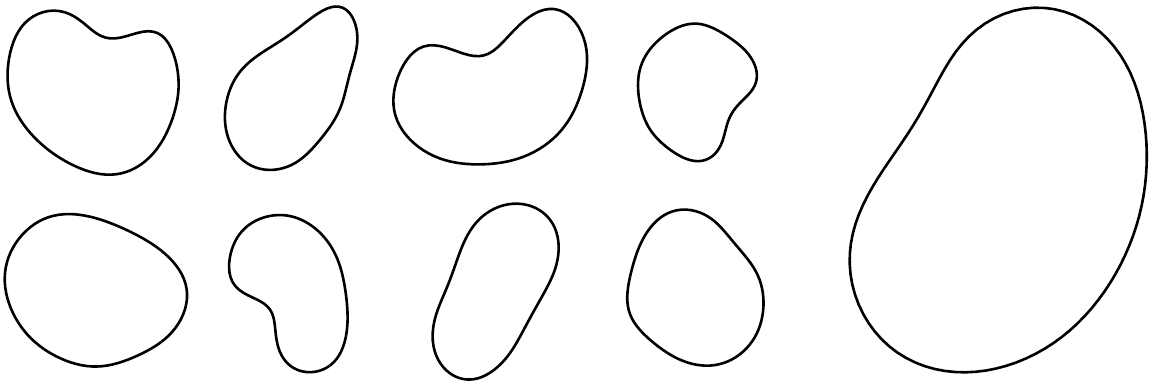}
\hspace{0.3cm}}}}
&
\includegraphics[width=0.21\textwidth]{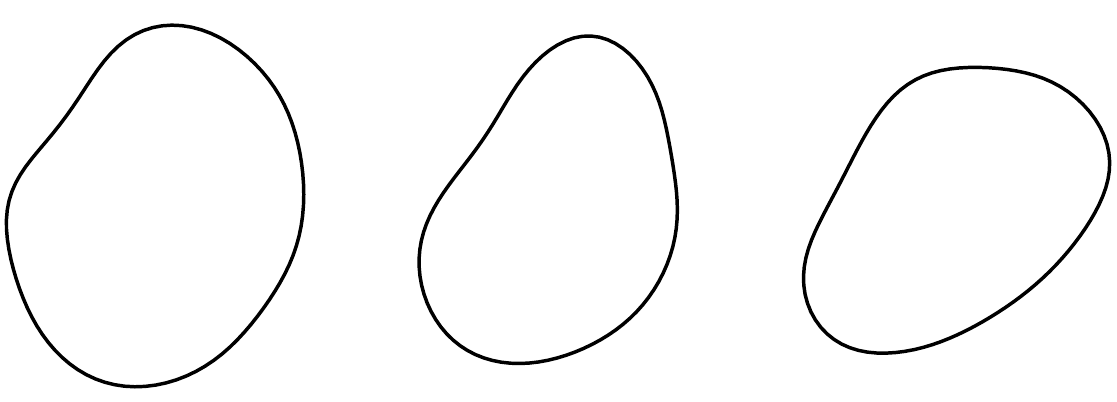}
\\
&
\includegraphics[width=0.21\textwidth]{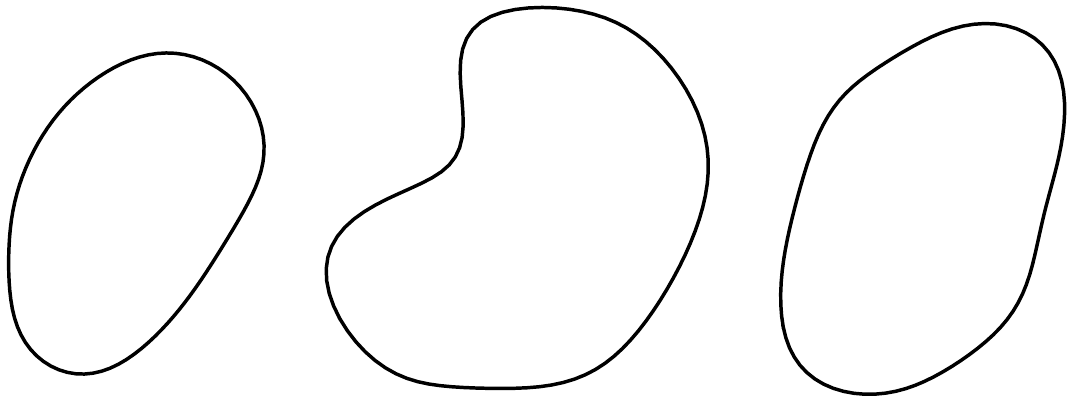}
\end{tabular}
\caption{Eight boundaries of HeLa cell nuclei, the Karcher mean of all nuclei (enlarged), and six randomly sampled cells using a Gaussian distribution in normal coordinates with the same covariance as the data.}
\label{hela_fig_1}
\end{figure}

To extract the boundary of the nucleus, we apply a thresholding method \cite{Otsu1979} to obtain a binary image and then fit -- using least squares -- a spline with $N_\th=12$ and $n_\th=4$ to the longest 4-connected component of the thresholded image. Then we reparametrize the boundary to approximately constant speed. The remaining degree of freedom is the starting point of the parametrization; when computing minimizing geodesics, we minimize over this parameter as well. The shapes of cell nuclei are thus represented by curves modulo translations, rotations, and constant shifts of the parametrisation, i.e. the two curves $c$ and $e^{i\be}c(\cdot - \al) + \la$ are considered equivalent. Examples of the extracted curves are depicted in Fig.~\ref{hela_fig_1}.

The parameters $a_0$, $a_1$, and $a_2$ of the Riemannian metric are chosen as follows: we compute the average $L^2$-, $H^1$- and $H^2$-contributions $\overline{E}_{L^2}$, $\overline{E}_{H^1}$, $\overline{E}_{H^2}$ to the energy of linear paths between each pair of curves in the dataset. As the $L^2$- and $H^2$-contributions scale differently, we rescale all curves such that $\overline{E}_{L^2} = \overline{E}_{H^2}$. Then we choose constants $a_0$, $a_1$, and $a_2$ such that
\[
a_0 \overline{E}_{L^2} : a_1 \overline{E}_{H^1} : a_2 \overline{E}_{H^2} = 3 : 1 : 6\;\text{ and }\;
\overline{E} = a_0 \overline{E}_{L^2} + a_1 \overline{E}_{H_1} +
a_2 \overline{E}_{H^2} = 100\enspace,
\]
resulting in an average length $\overline{\ell}_c = 12.45$ and $a_0 = 3.36$, $a_1 = 2.20$, and $a_2 = 6.73$.

\begin{figure}
\centering
\includegraphics[width=0.75\textwidth]{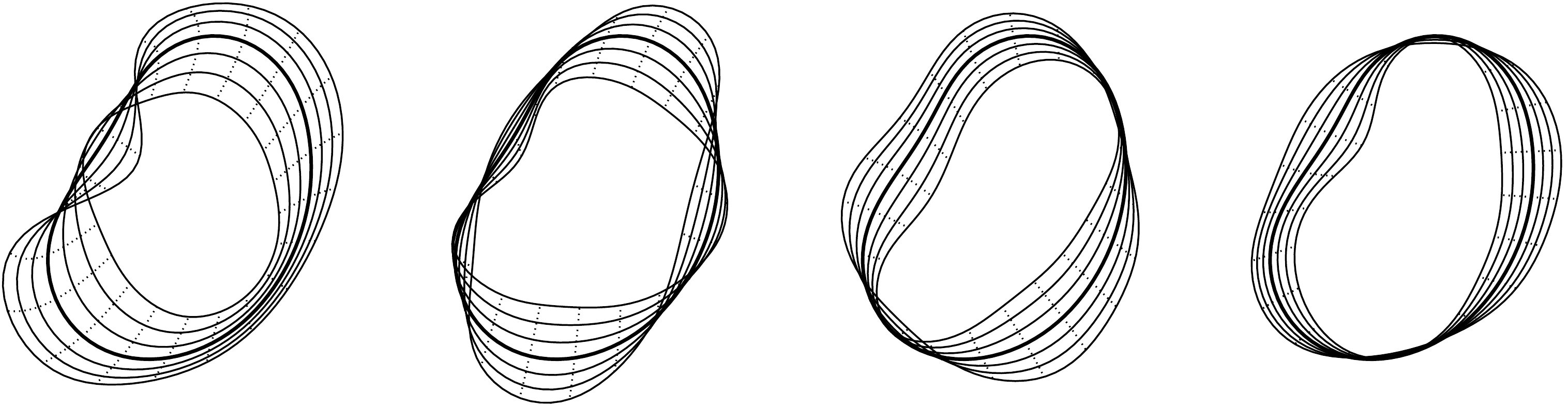}
\caption{Geodesics from the mean in the first four principal directions. The curves show the geodesic at times $-3, -2, \dots, 2, 3$; the bold curve is the mean. One can see four different characteristic deformations of the cell: bending, stretching along the long axis, stretching along the short axis, and a combination of stretching with partial bending.}
\label{hela_fig_2}
\end{figure}

The average shape of the nucleus can be captured by the Karcher mean $\overline{c}$. To compute the Karcher mean of the 87 nuclei, we use a conjugate gradient method on the Riemannian manifold of curves, as implemented in the Manopt library \cite{Manopt2014}, to solve the minimization problem \eqref{eq: KarcherEnergy}. We obtain convergence in $28$ steps; the final value of the objective function \eqref{eq: KarcherEnergy} is $F(\overline{c}) = 10.55$ and the norm of the gradient is $\| \on{grad}^G F(\overline{c})\|_{\overline c} < 10^{-3}$. The mean shape can be seen in Fig.~\ref{hela_fig_1}.

\begin{figure}
\centering
\includegraphics[width=0.6\textwidth]{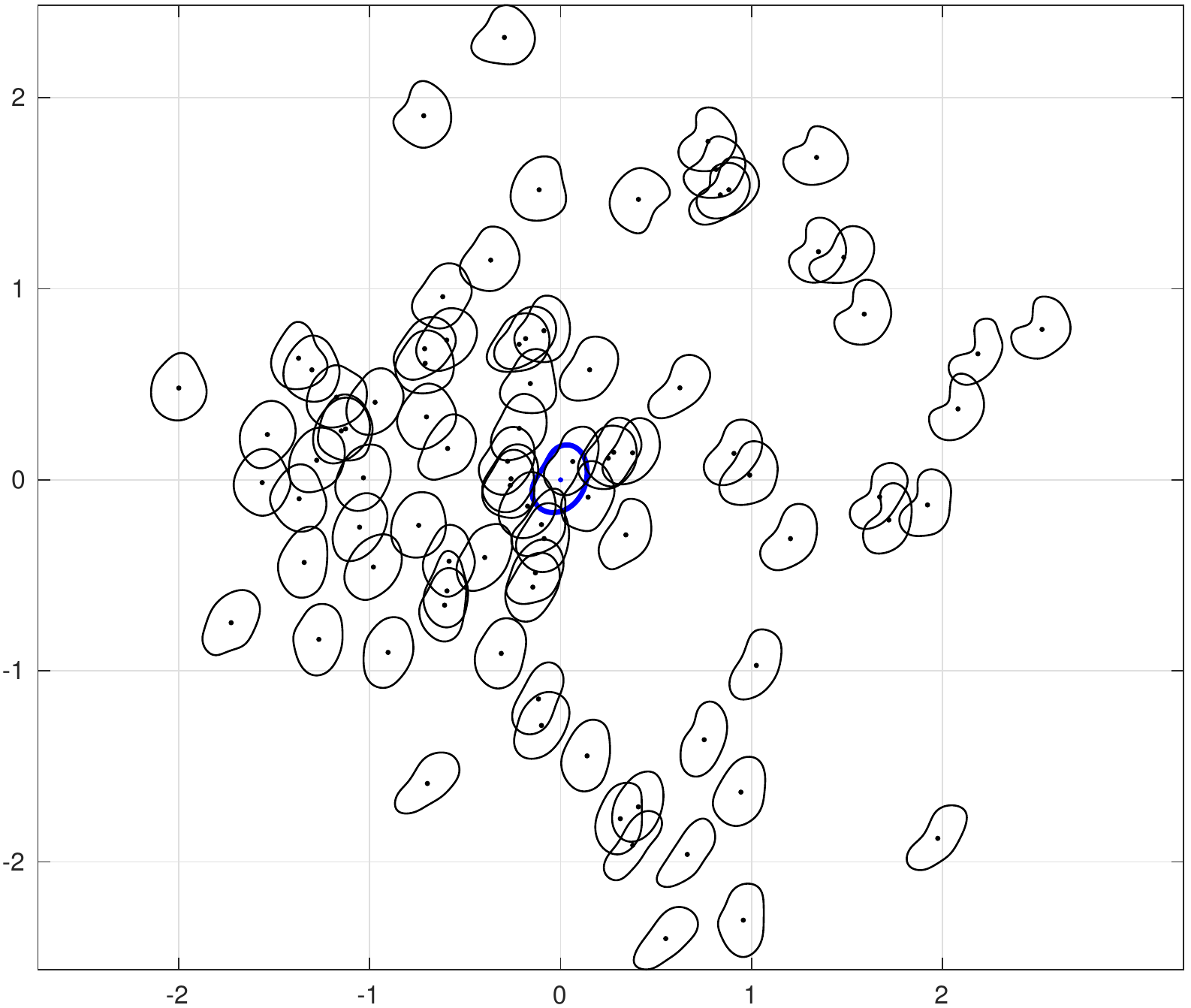}
\caption{All 87 cell nuclei projected to the plane in the tangent space which is spanned by the first two principal directions. The mean (in blue) is situated at the origin. The units on the coordinate axes are standard deviations. We see that the first coordinate is related to the bending of the nucleus, while the second coordinate is related to its elongation.}
\label{hela_fig_5}
\end{figure}

Having computed the mean $\overline{c}$, we represent each nuclear shape $c_j$ by the initial velocity $v_j = \on{Log}_{\overline c}c_j$ of the minimal geodesic connecting $\overline c$ and $c_j$. We perform principal component analysis with respect to the inner product $G_{\overline{c}}$ on the set of initial velocities $\{v_j :j = 1,\dots,87\}$. The first four eigenvalues are $4.10, 2.39, 1.68$ and $1.00$, and they explain $38.04\%,60.21\%,75.78\%$, and $85.07\%$ of the total variance. Geodesics from the mean in the directions of the first four principal directions can be seen in Fig.~\ref{hela_fig_2}. Fig.~\ref{hela_fig_5} shows the data projected to the subspace spanned by the first two eigenvectors. 

Finally, we sample from a normal distribution with the same covariance matrix as the data and project the sampled velocities $\hat v$ back to the space of curves using the exponential map $\hat c = \on{Exp}_{\overline{c}} \hat v$; some examples can be seen in Fig.~\ref{hela_fig_1}.

\subsection{Traces of Cardiac Images}\label{cardiac}

In our second application we study curves that are obtained from images of the cardiac cycle. More precisely, we consider a sequence of 30 cardiac images, taken at equispaced time points along the cardiac cycle. Each image is projected to a barycentric subspace of dimension two, yielding a closed curve in the two-dimensional space of barycentric coordinates. After normalizing the coordinates \cite[Sect.~3]{Pennec2015gsi} we obtain a closed, plane curve -- with the curve parameter representing time -- to which we can apply the methods presented in Sect.~\ref{math}. Details regarding the acquisition and projection of the images can be found in \cite{tobongomez,Mcleod}; barycentric subspaces on manifolds are described in \cite{Pennec2015gsi}.

\begin{figure}
	\centering
	\includegraphics[width=.15\textwidth]{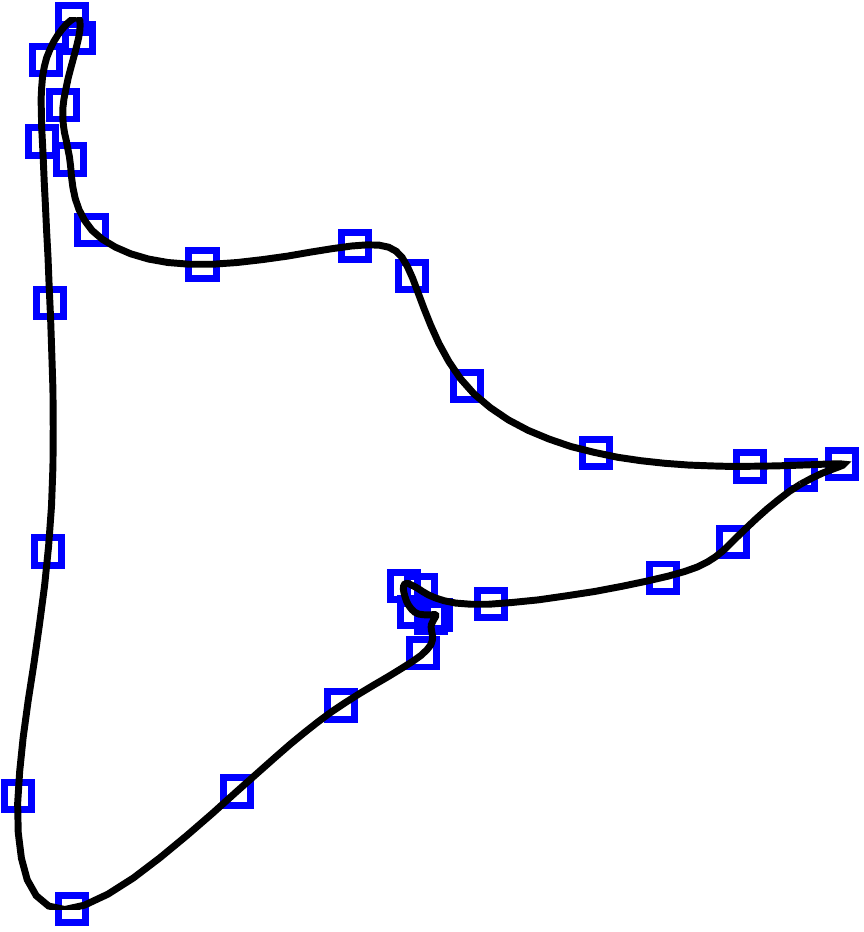}
	\hspace{1cm}
	\includegraphics[width=.15\textwidth]{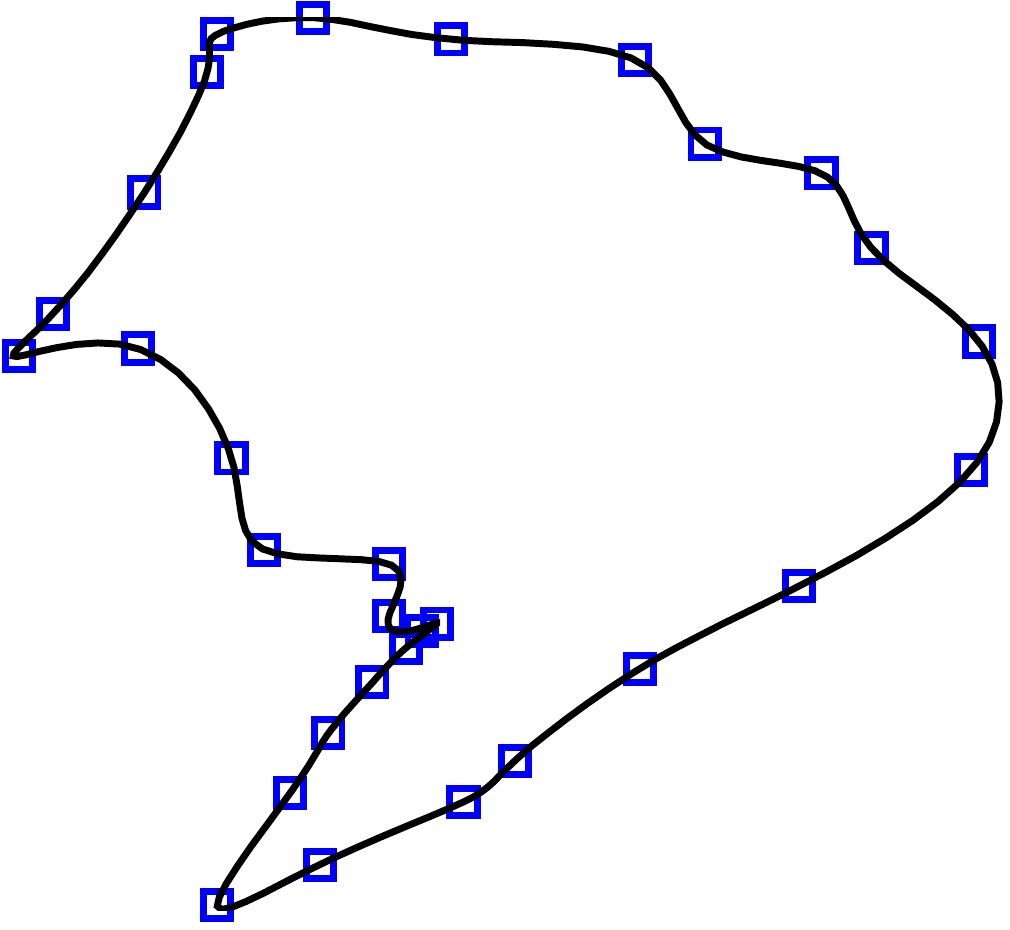}
	\hspace{1cm}	
	\includegraphics[width=.15\textwidth]{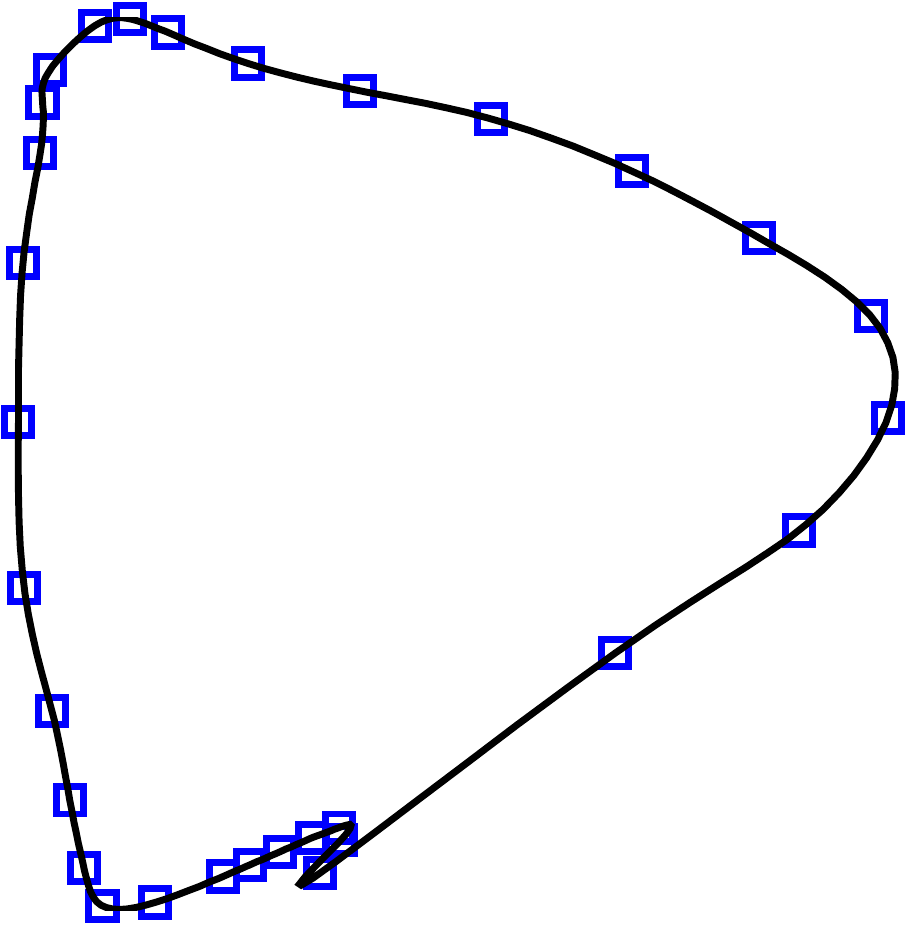}
	\caption{Projections to a two-dimensional barycentric subspace of 30 images from the cardiac cycles of three patients. Cubic splines interpolation of degree $n_\th=3$ with $N_\th=30$ control points is used.}
	\label{fig:splinefit}
\end{figure}

The data consists of 10 cardiac cycles of patients with Tetralogy of Fallot and 9 patients from a control group. Each cardiac cycle is originally represented by three-dimensional homogeneous coordinates $x_1:x_2:x_3$, sampled at 30 time points. We project   the homogeneous coordinates onto the plane $x_1+x_2+x_3=1$ and choose a two-dimensional coordinate system for this plane. Then we use spline interpolation with degree $n_\th=3$ and $N_\th=30$ control points to reconstruct the planar curves from the data points; see Fig.~\ref{fig:splinefit}.

\begin{figure}[h]
\centering
\subfloat[]{
	\includegraphics[width=.4\textwidth]{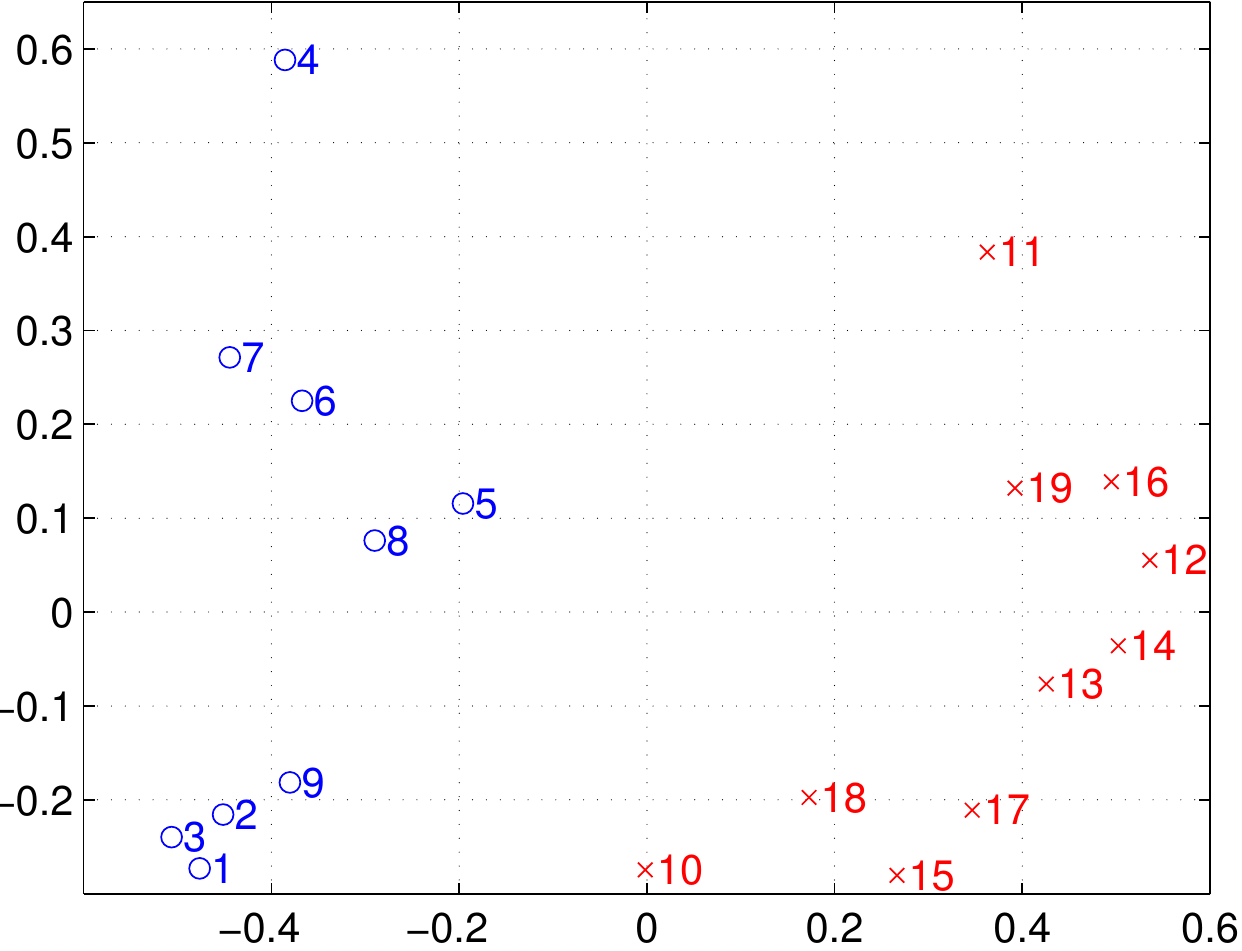}
}
\subfloat[]{
	\includegraphics[width=.4\textwidth]{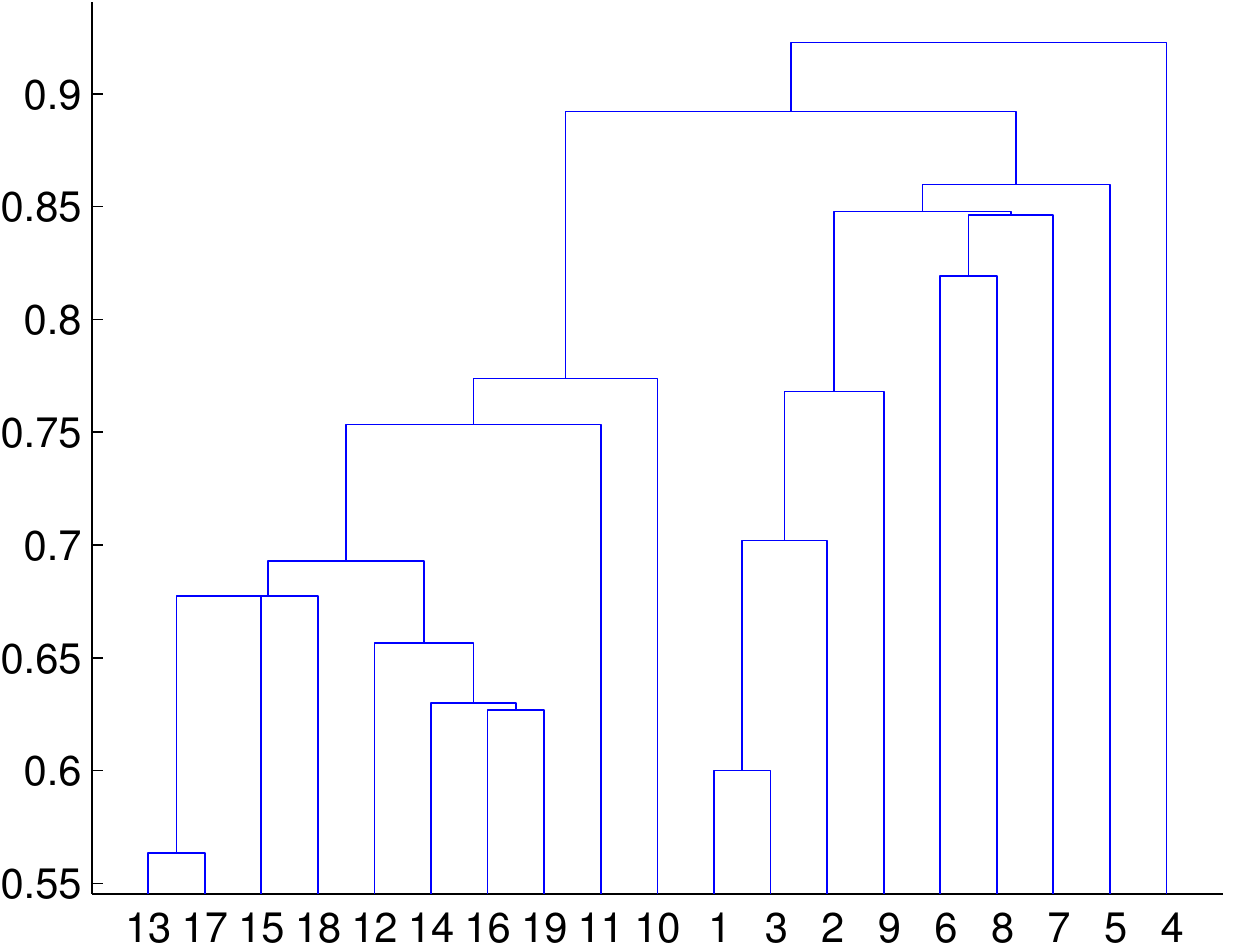}
}
\caption{(a) Two dimensional representation of the data using multi dimensional scaling of the pairwise distance matrix. (b) A dendrogram of clusters computed from the pairwise distance matrix using the single linkage criterion. Healthy patients are labelled 1--9 and diseased ones 10--19.}
\label{distanceMatrix}
\end{figure}

The parameters $a_0$, $a_1$, and $a_2$ in the metric are chosen similarly to Sect.~\ref{hela}; however, the scale of the curves is not changed and  
we use equal weighting between the $L^2$-, $H^1$- and $H^2$-parts of the average energy for linear paths. 
This leads to parameters $a_0 = 1$, $a_1 = 0.1$, and $a_2 = 10^{-9}$. To see if the metric structure derived from the Sobolev metric enables us to distinguish between diseased patients and the control group, we compute all 171  pairwise distances between the 19 curves; this takes about 15 minutes on a 2 GHz single core processor. Multi-dimensional scaling of the distance matrix shows that the metric separates healthy and diseased patients quite well (Fig.~\ref{distanceMatrix}a). Indeed, a cluster analysis based on the distance matrix recovers exactly -- with exception of one outlier (patient 4) -- the subgroups of healthy and diseased patients (Fig.~\ref{distanceMatrix}b). 

\begin{figure}
\centering
\includegraphics[width=.2\textwidth]{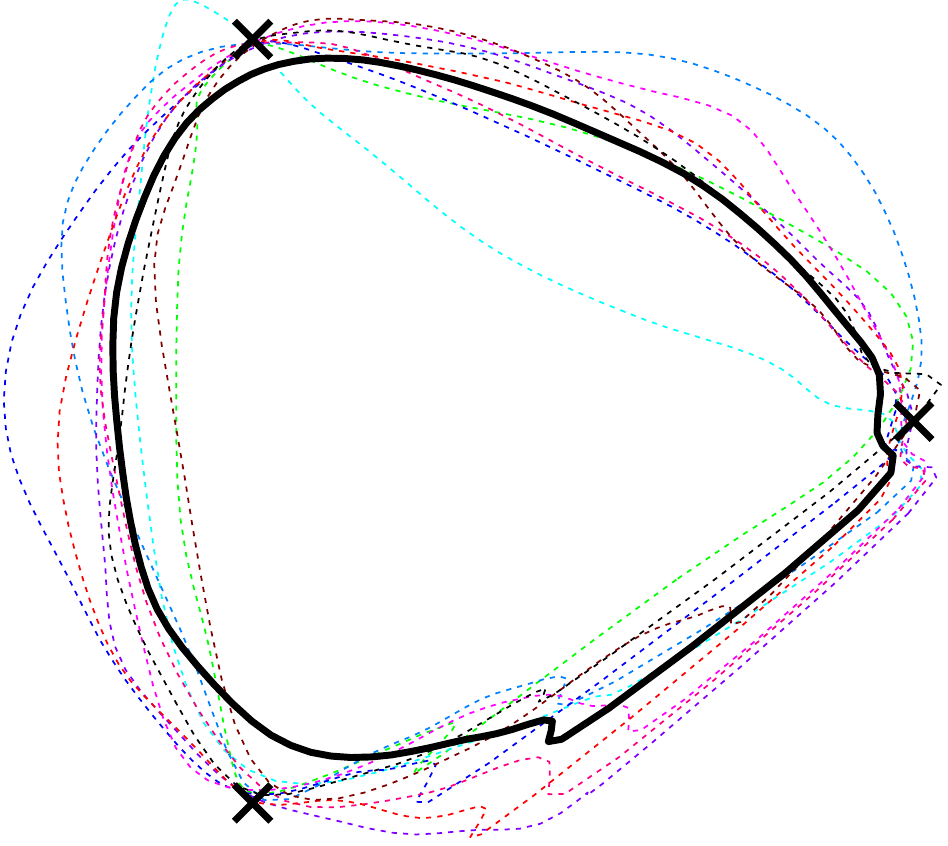}
\hspace{1cm}
\includegraphics[width=.2\textwidth]{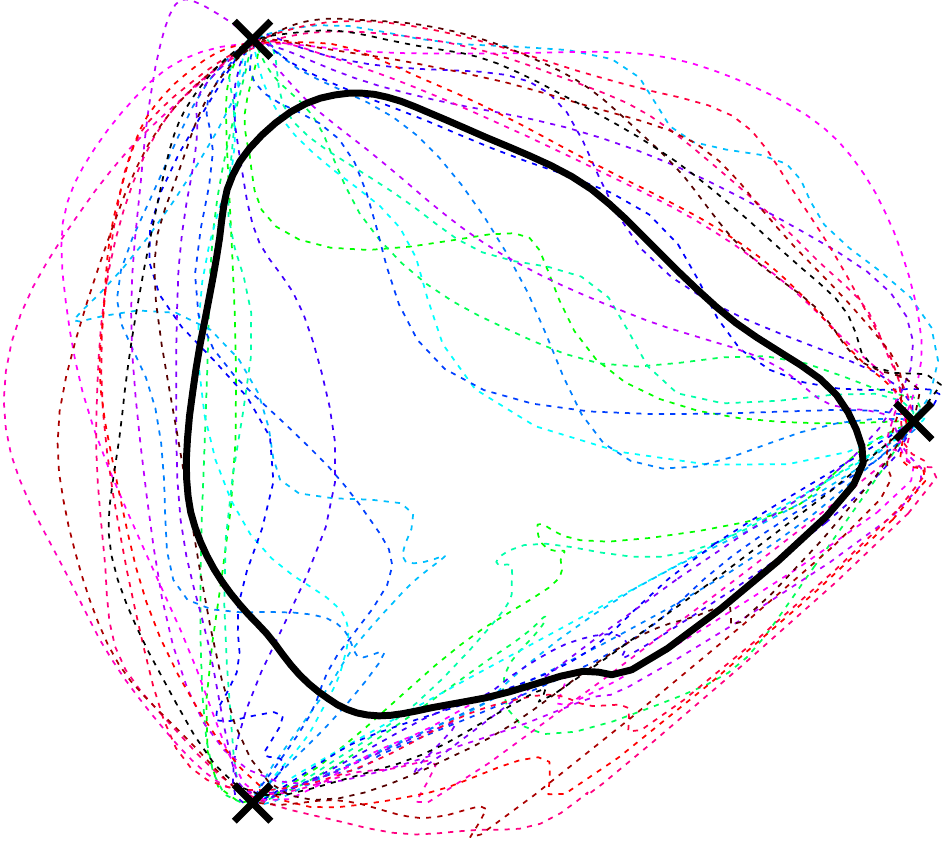}
\hspace{1cm}
\includegraphics[width=.2\textwidth]{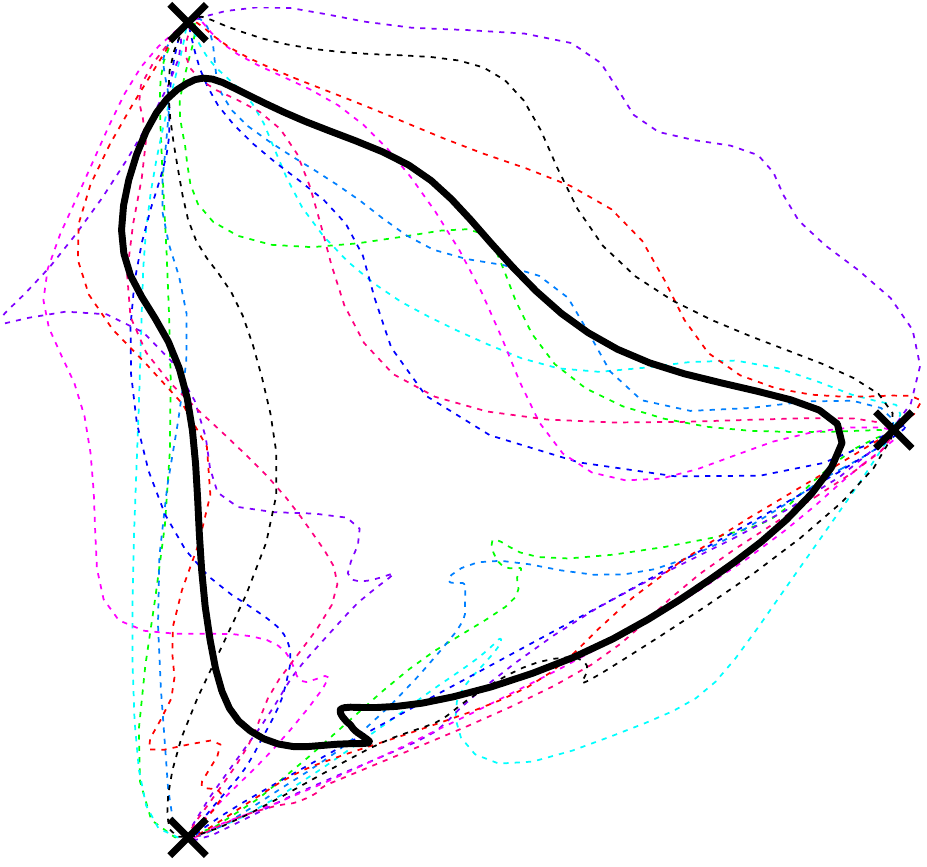}
\includegraphics[width=.75\textwidth]{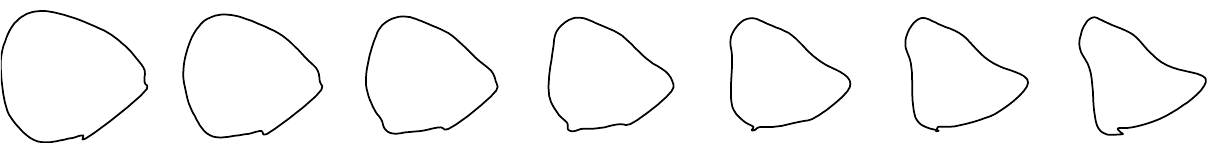}
\caption{First row: Karcher means of pathological cardiac cycles (left), all cycles (middle), and healthy cycles (right).
Second row: geodesic connecting the Karcher mean of pathological cycles to the Karcher mean of healthy cycles. The crosses denote the position of images, with respect to whom the barycentric projection was computed.}
\label{cardiacKarcher}
\end{figure}

The Karcher means of the healthy and diseased subgroups as well as of the entire population are depicted in Fig.~\ref{cardiacKarcher}. The mean was computed using a gradient descent method as described in Sect.~\ref{karcher} with a threshold of $10^{-4}$ for the norm of the gradient. 
The average distance from the mean for the diseased group is $0.6853$ with a variance of $0.0149$, and for the control group the distance is $0.7708$ with a variance of $0.0083$. 

\begin{figure}
\centering
\subfloat[]{
	\includegraphics[width=.31\textwidth]{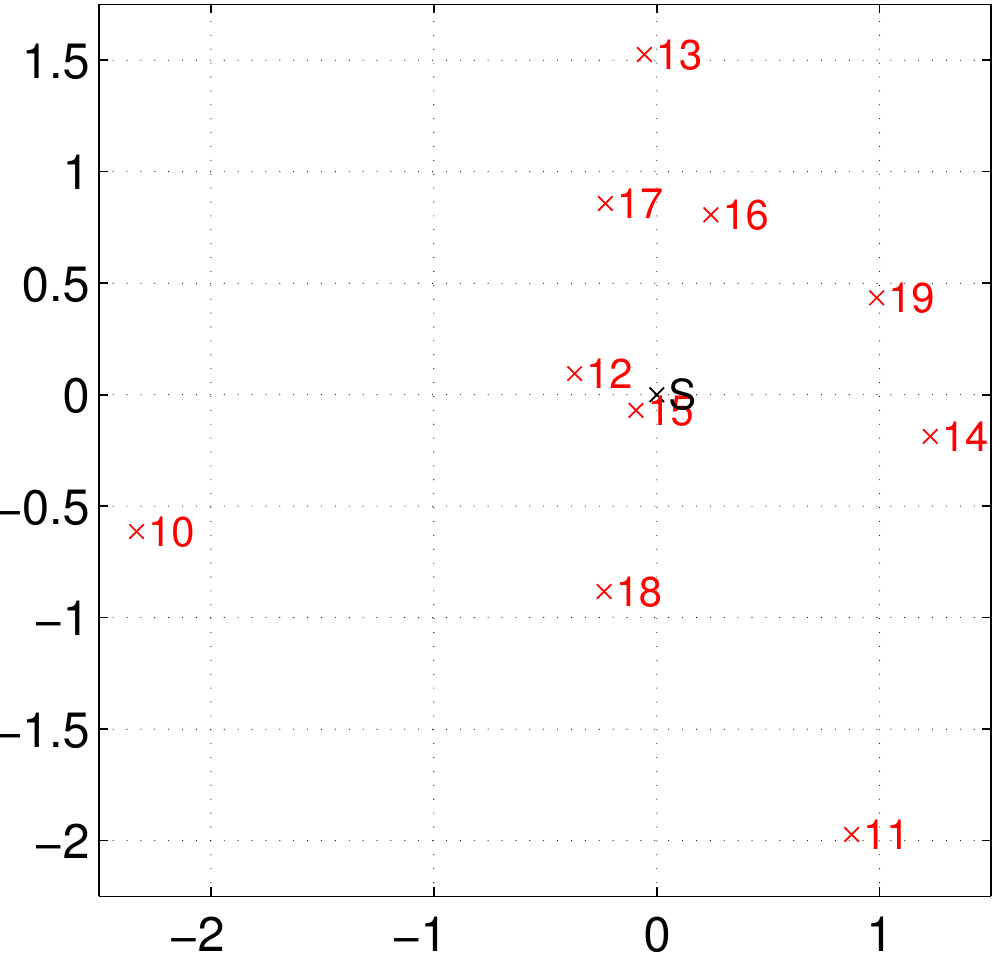}
}
\subfloat[]{
	\includegraphics[width=.31\textwidth]{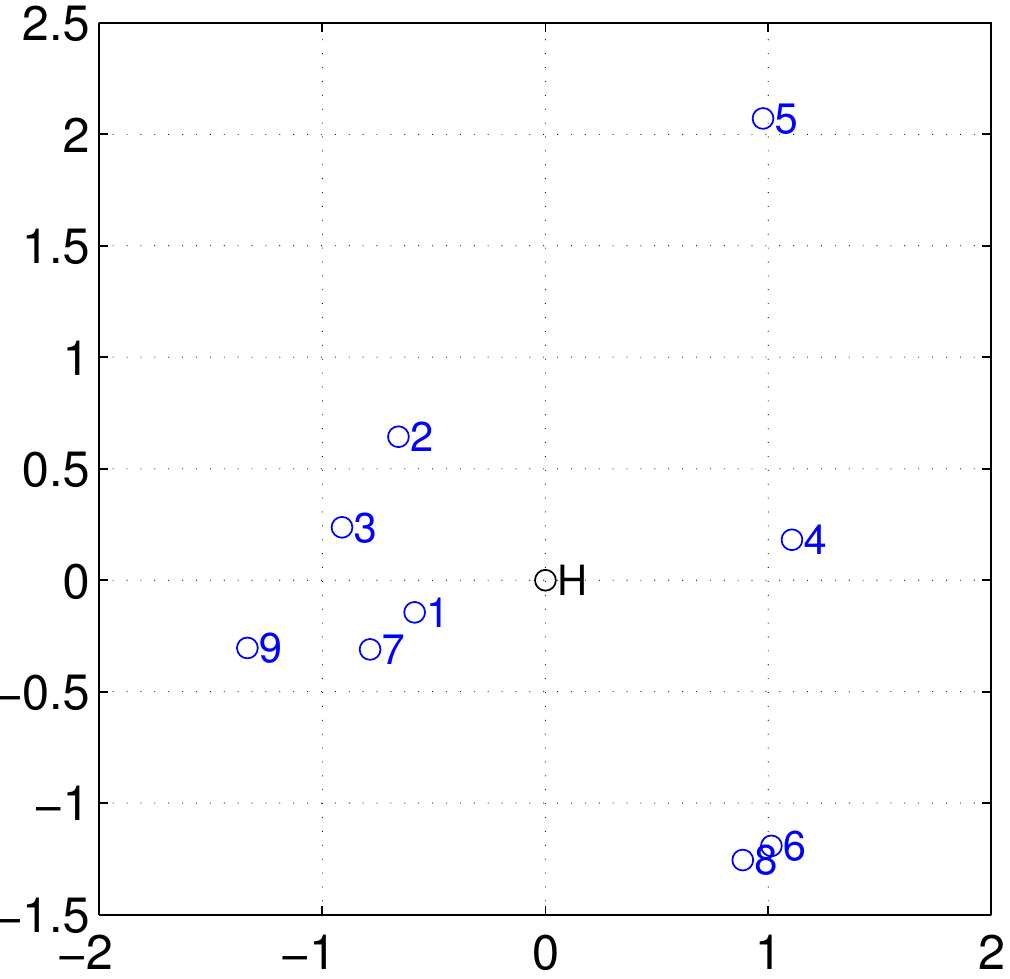}
}
\subfloat[]{
	\includegraphics[width=.31\textwidth]{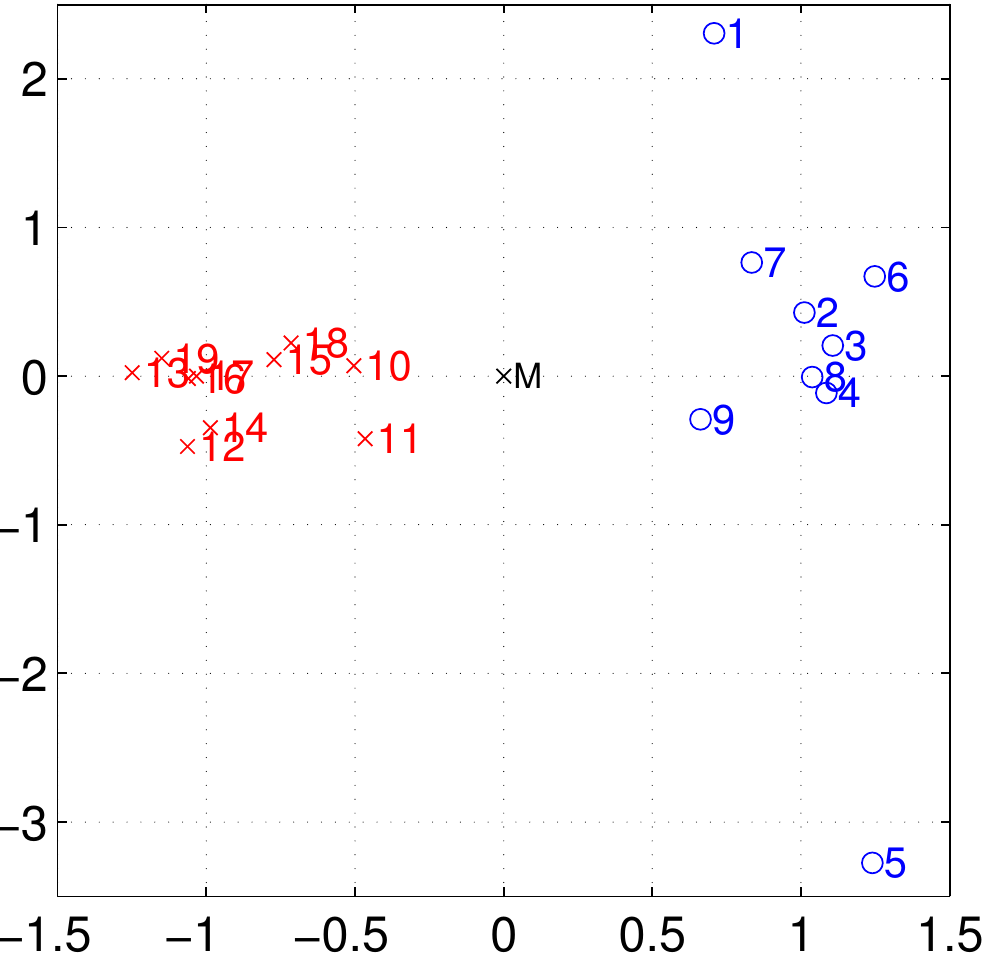} 
}
\caption{(a) Initial velocities of minimizing geodesics projected to the subspace spanned by first two principal components for the diseased group. (b) The same picture for the control group.
(c) The same picture for the whole population.}\label{cardiacKarcher2}
\end{figure}

To investigate the variability of the observed data, we performe principal component analysis on the initial velocities of the minimal geodesics connecting curves to the respective means (c.f. Sect.~\ref{hela}). Fig.~\ref{cardiacKarcher2} shows the initial velocities projected to the subspace spanned by the first two principal directions. Within the healthy and sick subgroups, less then $30\%$ of the principal components are needed to explain $90\%$ of the shape variation. If, in contrast, principal components are analyzed for the entire dataset based on the global Karcher mean, then $35\%$ of the principal components are needed to explain $90\%$ of the shape variation.

\section{Conclusions}
In this article we  numerically solved the initial and boundary value problems for geodesics on the space of regular curves  under second order Sobolev metrics. 
We analyzed two medical datasets using our approach.  In future work we plan to prove rigorous convergence results for our discretizations,
further investigate the impact of the constants in the metric, and treat unparametrized curves.

\subsection*{Acknowledgments}
We would like to thank Xavier Pennec and Marc-Michel Roh\'e for providing us the cardiac image data and Hermann Schichl for his invaluable help with AMPL.

\bibliographystyle{plain}

\begin{thebibliography}{10}

\bibitem{Bauer2011a}
M.~Bauer and M.~Bruveris.
\newblock {A new Riemannian setting for surface registration}.
\newblock In {\em {3nd MICCAI Workshop on Mathematical Foundations of
  Computational Anatomy}}, pages 182--194, 2011.

\bibitem{Bauer2014}
Martin Bauer, Martins Bruveris, and Peter~W. Michor.
\newblock Overview of the geometries of shape spaces and diffeomorphism groups.
\newblock {\em J. Math. Imaging Vis.}, 50:60--97, 2014.

\bibitem{Bauer2014c}
Martin Bauer, Martins Bruveris, and Peter~W. Michor.
\newblock {$R$}-transforms for {S}obolev {$H^2$}-metrics on spaces of plane
  curves.
\newblock {\em Geom. Imaging Comput.}, 1(1):1--56, 2014.

\bibitem{Bauer2011b}
Martin Bauer, Philipp Harms, and Peter~W. Michor.
\newblock {Sobolev metrics on shape space of surfaces}.
\newblock {\em J. Geom. Mech.}, 3(4):389--438, 2011.

\bibitem{Boland2001}
Michael~V. Boland and Robert~F. Murphy.
\newblock A neural network classifier capable of recognizing the patterns of
  all major subcellular structures in fluorescence microscope images of {HeLa}
  cells.
\newblock {\em Bioinformatics}, 17(12):1213--1223, 2001.

\bibitem{Manopt2014}
Nicolas Boumal, Bamdev Mishra, P.-A. Absil, and Rodolphe Sepulchre.
\newblock {M}anopt, a {M}atlab toolbox for optimization on manifolds.
\newblock {\em J. Mach. Learn. Res.}, 15:1455--1459, 2014.

\bibitem{Bruveris2014b_preprint}
Martins Bruveris.
\newblock Completeness properties of {S}obolev metrics on the space of curves.
\newblock {\em J. Geom. Mech.}, 7(2), 2015.

\bibitem{Bruveris2014}
Martins Bruveris, Peter~W. Michor, and David Mumford.
\newblock Geodesic completeness for {S}obolev metrics on the space of immersed
  plane curves.
\newblock {\em Forum Math. Sigma}, 2:e19, 2014.

\bibitem{Dryden1998}
I.~L. Dryden and K.~V. Mardia.
\newblock {\em Statistical Shape Analysis}.
\newblock Wiley Series in Probability and Statistics. John Wiley \& Sons, Ltd.,
  Chichester, 1998.

\bibitem{fourer2002ampl}
Rober Fourer, D~Gay, and Brian~W Kernighan.
\newblock The {AMPL} book, 2002.

\bibitem{Krim2006}
Hamid Krim and Anthony Yezzi, Jr., editors.
\newblock {\em Statistics and Analysis of Shapes}.
\newblock Modeling and Simulation in Science, Engineering and Technology.
  Birkh\"auser Boston, 2006.

\bibitem{Mcleod}
Kristin McLeod, Maxime Sermesant, Philipp Beerbaum, and Xavier Pennec.
\newblock {S}patio-temporal tensor decomposition of a polyaffine motion model
  for a better analysis of pathological left ventricular dynamics.
\newblock {\em IEEE Trans. Med. Imaging}, 2015.

\bibitem{Michor2007}
P.~W. Michor and D.~Mumford.
\newblock {An overview of the {R}iemannian metrics on spaces of curves using
  the {H}amiltonian approach}.
\newblock {\em Appl. Comput. Harmon. Anal.}, 23(1):74--113, 2007.

\bibitem{Mio2007}
Washington Mio, Anuj Srivastava, and Shantanu Joshi.
\newblock On shape of plane elastic curves.
\newblock {\em Int. J. Comput. Vision}, 73(3):307--324, July 2007.

\bibitem{Vialard2014_preprint}
G~Nardi, G~Peyr{\'e}, and F.-X. Vialard.
\newblock Geodesics on shape spaces with bounded variation and {S}obolev
  metrics.
\newblock {\em http://arxiv.org/abs/1402.6504}, 2014.

\bibitem{Otsu1979}
Nobuyuki Otsu.
\newblock A threshold selection method from gray-level histograms.
\newblock {\em IEEE T. Syst. Man Cyb.}, 9(1):62--66, 1979.

\bibitem{Pennec2006b}
Xavier Pennec.
\newblock Intrinsic statistics on {R}iemannian manifolds: basic tools for
  geometric measurements.
\newblock {\em J. Math. Imaging Vision}, 25(1):127--154, 2006.

\bibitem{Pennec2015gsi}
Xavier Pennec.
\newblock Barycentric subspaces and affine spans in manifolds, 2015.
\newblock To appear in the proceeding of Geometric Science of Information,
  2015.

\bibitem{Rohde2008}
Gustavo~K. Rohde, Alexandre J.~S. Ribeiro, Kris~N. Dahl, and Robert~F. Murphy.
\newblock Deformation-based nuclear morphometry: capturing nuclear shape
  variation in {HeLa} cells.
\newblock {\em Cytometry Part A}, 73A(4):341--350, 2008.

\bibitem{Rohde2008b}
Gustavo~K. Rohde, Wei Wang, Tao Peng, and Robert~F. Murphy.
\newblock Deformation-based nonlinear dimension reduction: applications to
  nuclear morphometry.
\newblock In {\em 5th IEEE Int. Symposium on Biomedical Imaging: From Nano to
  Macro}, pages 500--503, 2008.

\bibitem{Rumpf2014}
Martin Rumpf and Benedikt Wirth.
\newblock {Variational time discretization of geodesic calculus}.
\newblock {\em IMA Journal of Numerical Analysis}, 2014.

\bibitem{Schumaker2007}
Larry~L. Schumaker.
\newblock {\em Spline Functions: Basic Theory}.
\newblock Cambridge Mathematical Library. Cambridge University Press,
  Cambridge, third edition, 2007.

\bibitem{Jermyn2011}
Anuj Srivastava, Eric Klassen, Shantanu~H. Joshi, and Ian~H. Jermyn.
\newblock Shape analysis of elastic curves in {E}uclidean spaces.
\newblock {\em IEEE T. Pattern Anal.}, 33(7):1415--1428, 2011.

\bibitem{tobongomez}
Catalina Tobon-Gomez, Mathieu De~Craene, Kristin Mcleod, Lennart Tautz, Wenzhe
  Shi, Anja Hennemuth, Adityo Prakosa, Hengui Wang, Gerald Carr-White, Sergio
  Kapetanakis, Albert Lutz, Vernon Rasche, Tobias Schaeffter, Constantin
  Butakoff, Oskar Friman, Tommaso Mansi, Maxime Sermesant, Xiahai Zhuang,
  S{\'e}bastien Ourselin, Hans~Otto Peitgen, Xavier Pennec, Reza Razavi, Daniel
  Rueckert, Alejandro~F. Frangi, and Kawal Rhode.
\newblock {Benchmarking framework for myocardial tracking and deformation
  algorithms: an open access database}.
\newblock {\em {Medical Image Analysis}}, 17(6):632--648, 2013.

\bibitem{wachter2006implementation}
Andreas W{\"a}chter and Lorenz~T Biegler.
\newblock On the implementation of an interior-point filter line-search
  algorithm for large-scale nonlinear programming.
\newblock {\em Math. Program.}, 106(1):25--57, 2006.

\bibitem{Michor2008a}
L.~Younes, P.~W. Michor, J.~Shah, and D.~Mumford.
\newblock {A metric on shape space with explicit geodesics}.
\newblock {\em Atti Accad. Naz. Lincei Cl. Sci. Fis. Mat. Natur. Rend. Lincei
  (9) Mat. Appl.}, 19(1):25--57, 2008.

\bibitem{Younes2012}
Laurent Younes.
\newblock Spaces and manifolds of shapes in computer vision: an overview.
\newblock {\em Image Vision Comput.}, 30(6):389--397, 2012.

\end{thebibliography}

\end{document}